\documentclass[12pt,a4]{article}
\usepackage{amsmath}
\usepackage{amsfonts,amssymb}
\usepackage[]{graphicx}
\usepackage[cp1251]{inputenc}
\usepackage[english, russian]{babel}
\usepackage[margin=40mm]{geometry}
\begin{document}

\begin{flushleft}
 УДК 517.986.7:517.984.3
\end{flushleft}

\begin{center}
{\bf О НЕКОТОРЫХ СВОЙСТВАХ МНОГОМЕРНОГО ФУНКЦИОНАЛЬНОГО ИСЧИСЛЕНИЯ
БОХНЕРА-ФИЛЛИПСА\\
  А. Р. Миротин\\}
amirotin@yandex.ru
\end{center}

\hspace{5mm}

{\bf  1 Введение и предварительные сведения}

    Основы многомерного
исчисления Бохнера-Филлипса были заложены автором в [1--3],
\cite{M1}--\cite{M5}
 (относительно одномерного исчисления см. [4, 5], а также
[6--10]; отметим, что оно находит важные применения в теории
вероятностей, которые связаны с тем, что функции Бернштейна
возникают в теории безгранично делимых распределений [11,
глава~XIII, \S\ 7, теорема 2], а также теории случайных процессов
[12, глава~ XXIII, \S\ 5, теорема 23.115.2 и формула (23.15.13);
13, с. 1339, раздел 7]).
     В данной статье мы
определяем совместный спектр набора замкнутых коммутирующих
операторов в банаховом пространстве и устанавливаем теоремы об
отображении этого спектра и его частей для многомерного исчисления
Бохнера-Филлипса (теоремы об отображении  спектра для одномерного
исчисления см. в [9]), даем условие голоморфности полугрупп,
порождаемых операторами, возникающими в исчислении,  а также
доказываем полезное неравенство моментов для этих операторов.
Напомним необходимые понятия и факты из [1] и [3].

{\bf 1. Определение} [1]. Будем говорить, что неположительная
функция $\psi\in C^\infty((-\infty;0)^n)$ есть (отрицательная)
\textit{функция Бернштейна $n$ переменных} (или {\it принадлежит
классу ${\cal T}_n$}), если все ее частные производные первого
порядка абсолютно монотонны (функция из
$C^{\infty}((-\infty;0)^n)$ называется {\it абсолютно монотонной,}
если она неотрицательна вместе со своими частными производными
всех порядков).

Ясно, что  ${\cal T}_n$ есть конус относительно поточечного
сложения функций и умножения на скаляр.

Известно [3], что каждая функция $\psi\in {\cal T}_n$ допускает
интегральное представление (здесь и ниже точкой мы обозначаем
скалярное произведение в $\Bbb{R}_+^n$; запись $s\to -0$ означает,
что $s_1\to -0, \dots, s_n\to -0$)

$$
\psi(s)=c_0+c_1\cdot s+\int\limits_{\Bbb{R}_+^n\setminus \{0\}}
(e^{s\cdot u}-1)d\mu (u)   (s\in(-\infty;0)^n),  \eqno(1.1)
$$

\noindent где $c_0=\psi(-0):=\lim\limits_{s\to -0}\psi(s)$, а
$c_1$ из $\Bbb{R}_+^n$ и положительная мера $\mu$ на
$\Bbb{R}_+^n\setminus \{0\}$ определяются по $\psi$ однозначно.
Кроме того, $\psi$  продолжается по формуле (1.1) до функции,
голоморфной в области $\{\mathrm{Re}s<0\}\subset \mathbb{C}^n$, и
непрерывной в замыкании этой области.

Примерами функций Бернштейна одного переменного  служат
$-(-s)^\alpha$  $(0<\alpha\leq 1),\ -\log(1-s),\
-\mathrm{arch}(1-s)$ (см. [9]). Большое число примеров приведено в
монографии [10], специально посвященной этому предмету.

Всюду далее, если не оговорено противное,  $T(u)=T_1(u_1)\dots
T_n(u_n)$ $(u\in \Bbb{R}_+^n)$ есть ограниченная
$n$-параметрическая $C_0$-полугруппа операторов  в комплексном
банаховом пространстве $X$,
 $T_1, \dots ,T_n$ --- попарно коммутирующие
однопараметрические $C_0$-полугруппы  в  $X$, удовлетворяющие
условиям $||T_j(t)||\leq M_j\quad(t\geq 0, M_j={\rm const}\geq 1,
\ j=1, \dots ,n)$. Через $A_j$ обозначим генератор полугруппы
$T_j$ с областью определения $D(A_j)$ (частный генератор
полугруппы $T$) и положим $A=(A_1,\dots ,A_n)$. Далее
коммутирование операторов $A_1,\dots, A_n$ означает коммутирование
соответствующих полугрупп. Через ${\rm Gen}(X)$ мы будем
обозначать множество всех генераторов равномерно ограниченных
$C_0$-полугрупп в $X$, а через $I$ -- единичный оператор в $X$.
Отметим, что  векторное подпространство $D(A):=\cap_{j=1}^nD(A_j)$
плотно в $X$ (см. [12, \S\,10.10] или [14, предложение 3.26]).

{\bf 2. Определение} [3]. Определим значение функции $\psi$ из
${\cal T}_n$ вида (1.1) на наборе $A=(A_1,\dots ,A_n)$ при $x\in
D(A)$ формулой
$$
\psi(A)x=c_0x+c_1\cdot Ax+\int\limits_{\Bbb{R}_+^n\setminus \{0\}}
(T(u)-I)xd\mu(u),     \eqno(1.2)
$$
где $c_1\cdot Ax:=\sum_{j=1}^nc_1^jA_jx$ (интеграл понимается в
смысле Бохнера).

Пусть $\psi\in{\cal T}_n, t\geq 0$. Тогда функция
$g_t(s):=e^{t\psi(s)}$ будет абсолютно монотонной на
$(-\infty;0)^n$, и $g_t(s)\leq 1.$ В силу многомерного варианта
теоремы Бернштейна-Уиддера (см., например, [15]) существует такая
единственная ограниченная положительная  мера $\nu_t$ на
$\Bbb{R}_+^n,$ что при $s\in (-\infty;0)^n$
$$
g_t(s)=\int\limits_{\Bbb{R}_+^n} e^{s\cdot u}d\nu_t(u)={\cal
L}\nu_t(s)\eqno(1.3)
$$
(здесь и ниже ${\cal L}\nu$ обозначает $n$-мерное преобразование
Лапласа меры  $\nu$).

{\bf 3. Определение.} Используя обозначения, введенные выше,
положим ($x\in X$)
$$
g_t(A)x=\int\limits_{\Bbb{R}_+^n} T(u)xd\nu_t(u)\eqno(1.4)
$$
(интеграл  Бохнера).

Очевидно, что $\|g_t(A)\|\leq (\prod_{j=1}^nM_j) e^{t\psi(-0)}\leq
\prod_{j=1}^nM_j.$ Поскольку $g_{t+r}(s)=g_t(s)g_r(s),$ то $\nu_t$
образуют сверточную полугруппу ограниченных  мер на $\Bbb{R}_+^n$.
Поэтому $g(A):t\mapsto g_t(A)$ есть  ограниченная полугруппа
операторов на $X$. В частности, $g(A)$ есть $C_0$-полугруппа. В
одномерном случае она называется {\it полугруппой, подчиненной
полугруппе $T$} (терминология восходит к теории вероятностей, см.
[11, \S\,X.7].

В [16] было замечено, что  основная теорема из [3] может быть
уточнена следующим образом.

{\bf 4. Теорема}. {\it Замыкание оператора $\psi(A)$ существует и
является генератором  полугруппы $g(A)$ класса $C_0$, определенной
формулой \rm{(1.4)}.}

В самом деле, в [3] было доказано, что оператор $\psi(A)$
замыкаем, и его расширением является генератор $G$
$C_0$-полугруппы $g(A)$. Поэтому замыкание
$\overline{\psi(A)}\subseteq G$. Покажем, что здесь имеет место
равенство. Поскольку операторы $g_t(A)$ коммутируют с $T_k(s)$ при
всех $k$, то, как легко проверить, $g_t(A):D(A_k)\to D(A_k)$ при
всех $k$, а потому и $g_t(A):D(A)\to D(A)$. Отсюда следует, что
$D(A)$ есть существенная область для генератора $G$ (см. [17,
следствие 3.1.7]). С другой стороны, $D(A)$ есть существенная
область для оператора $\overline{\psi(A)}$, причем сужения
операторов $\overline{\psi(A)}$ и $G$ на $D(A)$ совпадают. Поэтому
$\overline{\psi(A)}= G$,
 что и завершает доказательство.$\Box$

Предыдущая теорема мотивирует окончательный вариант определения
операторов $\psi(A)$.

{\bf 5. Определение } [3]. Под {\it значением функции} $\psi$ из
${\cal T}_n$ на наборе $A=(A_1,\dots ,A_n)$ коммутирующих
операторов из ${\rm Gen}(X)$ будем понимать генератор полугруппы
$g(A)$. Это значение мы далее обозначаем $\psi(A)$.  Возникающее
функциональное исчисление будем называть {\it многомерным
исчислением Бохнера-Филлипса}, или {\it ${\cal T}_n$-исчислением.}

Введенные выше обозначения и ограничения далее будут применяться
без дополнительных пояснений.

\bigskip
 {\bf 2 Отображение совместных спектров}

Уже для функции $\psi\in {\cal T}_1$ теорема об отображении
спектра в смысле справедливости равенства
$\psi(\sigma(A))=\sigma(\psi(A))$, где $A\in {\rm
Gen}(X),\sigma(A)$ -- спектр оператора $A$, может не выполняться,
в чем нас убеждает следующий пример.

{\bf 6. Пример}. Пусть $\psi(s)= e^s-1.$ В этом случае
представляющая мера $\mu$ равна $\varepsilon_1$, мере Дирака на
$\Bbb{R}$, сосредоточенной в единице. В силу определения 5 и
теоремы 4 тогда $\psi(A)=T(1)-I$. Но известно, что существуют
такие равномерно ограниченные полугруппы класса $C_0$, для которых
$\sigma(T(1))\ne e^{\sigma(A)}$ (см., например, [18, глава~1,
пример 9.6]; в этом примере $\sigma(A)=\emptyset$, а $T(t)=O$ при
$t\geq 1$). Последнее означает, что $\sigma(\psi(A))\ne
\psi(\sigma(A)).$ $\Box$

\bigskip
{\bf 7. Определение}. {\it Совместный точечный спектр}
$\sigma_p(A)$ набора $A$  операторов в $X$ есть множество тех
$\lambda\in \Bbb{C}^n$, для которых найдется вектор $x\in D(A),
x\ne 0$ (\textit{совместный собственный вектор набора }$A$),
удовлетворяющий равенствам $A_jx=\lambda_jx$ при всех
$j=1,\dots,n$.

\bigskip
{\bf 8. Определение}. {\it Совместный аппроксимативный спектр}
$\sigma_a(A)$ набора $A$  операторов в $X$ есть множество тех
$\lambda\in \Bbb{C}^n$, для которых найдется такая
последовательность векторов  $x_m\in D(A), \|x_m\|=1$
(\textit{совместный аппроксимативный собственный вектор набора
}$A$), что $\|A_jx_m-\lambda_jx_m\|\to 0\ (m\to\infty)$ при всех
$j=1,\dots,n$.

\bigskip
{\bf 9. Определение}. {\it Совместный остаточный спектр}
$\sigma_R(A)$ набора $A$  операторов в $X$ есть множество тех
$\lambda\in \Bbb{C}^n$, для которых векторное пространство
$\sum_{j=1}^n \mathrm{Im}(\lambda_j-A_j)$ не плотно в $X$.

\bigskip
{\bf 10. Определение}. {\it Совместный  спектр} $\sigma_J(A)$
набора $A$  операторов в $X$  определим равенством
$$
\sigma_J(A)=\sigma_a(A)\cup\sigma_R(A).
$$

Легко видеть, что в случае одного оператора эти определения
совпадают с соответствующими классическими, и $\sigma_p(A)$,
$\sigma_a(A)$, $\sigma_R(A)$ и $\sigma_J(A)$ содержатся в
$\sigma_p(A_1)\times\dots\times\sigma_p(A_n)$,
$\sigma_a(A_1)\times\dots\times\sigma_a(A_n)$,
$\sigma_R(A_1)\times\dots\times\sigma_R(A_n)$, и
$\sigma(A_1)\times\dots\times\sigma(A_n)$ соответственно.

Далее мы положим $A':=(A_1',\dots,A_n')$ (штрих здесь обозначает
сопряженный оператор).

\bigskip
{\bf 11. Лемма}. {\it Для любого набора $A$  операторов в $X$
справедливо равенство $\sigma_R(A)=\sigma_p(A')$.}

Доказательство. Заметим, что $\lambda\in\sigma_R(A)$ тогда и
только тогда, когда $f(\sum_{k=1}^n
\mathrm{Im}(\lambda_k-A_k))=\{0\}$, т. е. (поскольку $\sum_{k=1}^n
\mathrm{Im}(\lambda_k-A_k)\supseteq \mathrm{Im}(\lambda_j-A_j)$
при всех $j=1,\dots, n$) тогда и только тогда, когда
$f(\mathrm{Im}(\lambda_j-A_j))=\{0\}$ при всех $j=1,\dots, n$ для
некоторого ненулевого непрерывного линейного функционала $f$ в
$X$.  В силу определения сопряженного оператора это равносильно
тому, что $(\lambda_j-A'_j)f=0$ при всех $j=1,\dots, n. \Box$

\bigskip
{\bf  12. Теорема}. \textit{Для набора $A$ попарно коммутирующих
генераторов ограниченных $C_0$-полугрупп в банаховом пространстве
$X$ справедливы включения}

1) $\sigma_R(\psi(A))\supseteq \psi(\sigma_R (A))$;

2) $\sigma_p(\psi(A))\supseteq \psi(\sigma_p (A))$;

3) $\{\alpha\in\sigma_p(\psi(A)):\mbox{ найдутся собственный
вектор } x \mbox{ оператора } \psi(A),$
$\mbox{ отвечающий } \alpha, \mbox{ и совместный собственный }
\mbox{ вектор } f \mbox{ набора } A', \mbox{ такие},$   
$\mbox{
что } f(x)\ne 0 \}\subseteq \psi(\sigma_R(A)).$

\textit{Если к тому же $\sigma_a(A_j)\subseteq \{z\in
\mathbb{C}:\mathrm{Re}z<0\}$  при всех $j$,  или
$\partial\psi(-0)/\partial s_j\ne\infty$ при всех $j$, то }

4) $\sigma_a(\psi(A))\supseteq \psi(\sigma_a (A))$;

5) $\sigma(\psi(A))\supseteq \psi(\sigma_J (A))$.

\medskip
Доказательство.  1) При $\lambda\in\Bbb{C}^n, {\rm Re}\lambda\leq
0$ и $x\in D(A)$ имеем в силу (1.1) и (1.3)
$$
\psi(\lambda)x - \psi(A)x = c_1\cdot(\lambda I - A)x +
\int\limits_{\Bbb{R}_+^n} (e^{\lambda\cdot u}I -
T(u))xd\mu(u).\eqno(2.1)
$$

Для любых ограниченных операторов $B_j, C_j, j=1,\dots,n$ в $X$
справедливо тождество
$$
\prod_{j=1}^nB_j - \prod_{j=1}^nC_j = (B_1 - C_1)\prod_{j=2}^nB_j
+ C_1(B_2 - C_2)\prod_{j=3}^nB_j + \dots +
\prod_{j=1}^{n-1}C_j(B_n - C_n).
$$
С учетом коммутирования полугрупп $T_j$ получаем отсюда, что
$$
e^{\lambda\cdot u}I-T(u)=\sum_{j=1}^n(e^{\lambda_ju_j}I-T_j(u_j))
U_j(u),\eqno(2.2)
$$
\noindent где
$$
U_j(u)=\prod\limits_{1\leq l<j}T_l(u_l)\prod\limits_{j<k\leq
n}e^{\lambda_k u_k}
$$
--- ограниченные операторы в $X$, коммутирующие с $T_i
\ (i=1,\dots,n)$, $\|U_j(u)\|\leq M^{n-1}$.  Известно (см.,
например,  [14, формулы (8.1a), (8.1b)]), что  формулы
$$
V^\lambda_j(u_j)x=\int\limits_0^{u_j}e^{(u_j-s)\lambda_j}T_j(s)xds
$$
определяют ограниченные операторы в $X, {\rm
Im(}V^\lambda_j(u_j))\subseteq D(A_j),$ также коммутирующие со
всеми $T_i$, причем
$$
(e^{\lambda_ju_j}I-T_j(u_j))x=(\lambda_jI-A_j)V^\lambda_j(u_j)x
\mbox{ при } x\in X.
$$
 Подставляя это в (2.2), а результат -- в (2.1), имеем
при $x\in D(A)$
$$
(\psi(\lambda)I-\psi(A))x=\sum_{j=1}^nc_1^j(\lambda_jI-A_j)x+
\sum_{j=1}^n\int\limits_{\Bbb{R}_+^n}(\lambda_jI-A_j)V^\lambda_j(u_j)U_j(u)xd\mu(u).\eqno(2.3)
$$

Если $\lambda\in\sigma_R(A)$, то правая часть здесь принадлежит
замыканию в $X$ подпространства $\sum_{j=1}^n
\mathrm{Im}(\lambda_j-A_j)$, причем это замыкание не совпадает с
$X$. Следовательно, с учетом того, что $D(A)$ есть существенная
область оператора $\psi(A)$, получаем, что образ оператора
$\psi(\lambda)I-\psi(A)$ не всюду плотен, т. е. $\psi(\lambda)\in
\sigma_R(\psi(A))$.

2) Пусть $A_jx=\lambda_jx$ при некотором $x\in D(A), x\ne 0,
 j=1,\dots, n.$ Тогда $T(u)x=e^{\lambda\cdot u}x,$ поскольку в силу
[12, теорема 11.6.3], $T_j(\xi)x=e^{\lambda_j\xi}x.$ Формулы (1.3)
и (1.1) показывают теперь, что
$$
\psi(A)x=(c_0+c_1\cdot\lambda+\int\limits_{\Bbb{R}^n_+}(e^{\lambda\cdot
u}-1) d\mu(u))x=\psi(\lambda)x,
$$
 т. е. $\psi(\lambda)\in \sigma_p(\psi(A)),$ что и
требовалось.

3) Пусть $\alpha$ принадлежит левой части доказываемого включения,
вектор $x\in D(\psi(A)),\ x\ne 0$ таков, что $\psi(A)x=\alpha x$,
 и функционал $f\in X'$ таков, что $f(x)\ne 0$, и $A'f=\lambda f$
 для некоторого $\lambda\in \sigma_p(A')$. Покажем, что
 $\alpha=\psi(\lambda)$. С этой целью заметим, что в силу теоремы
 4 и определения 5 найдется такая последовательность $x_m\in D(A)$, что
 $x_m\to x$,  и $\psi(A)x_m\to \psi(A)x$. Следовательно, $\lim_mf(\psi(A)x_m)=
 \alpha f(x)$.  С учетом равенств $f(A_jy)=\lambda_jf(y)\ (j=1,\dots,n; y\in X)$
  последнее соотношение принимает вид
$$
\lim_{m\to\infty}
\left(c_0f(x_m)+(c_1\cdot\lambda)f(x_m)+\int\limits_{\Bbb{R}_+^n}
(f(T(u)x_m)-f(x_m))d\mu(u)\right)=\alpha f(x).\eqno(2.4)
$$

Для фиксированного $m$ рассмотрим функцию $\varphi(u)=f(T(u)x_m)\
(u\in \Bbb{R}_+^n)$. Так как   $T(u)x_m\in D(A)$, при любом $j$
имеем
$$
\frac{\partial\varphi(u)}{\partial s_j}=\lim_{h\to
+0}\frac{f(T_j(h)T(u)x_m)-f(T(u)x_m)}{h}=
$$
$$
 =f(A_jT(u)x_m)=
\lambda_jf(T(u)x_m)=\lambda_j\varphi(u).
$$
Отсюда $\varphi(u)=Ce^{\lambda\cdot u}$, где
$C=\varphi(0)=f(x_m)$. Значит, (2.4) можно записать как
$$
\lim_{m\to\infty}
\left(c_0+c_1\cdot\lambda+\int\limits_{\Bbb{R}_+^n}
(e^{\lambda\cdot u}-1)d\mu(u)\right)f(x_m)=\alpha f(x),
$$
откуда  $\alpha=\psi(\lambda)$. Применение леммы 11 завершает
доказательство.

 4) Заметим, что если $\mathrm{Re}\lambda_j<0$, то оператор
 $$
W^\lambda_jx=c_1^jx+\int\limits_{\Bbb{R}_+^n}V^\lambda_j(u_j)U_j(u)xd\mu(u)
$$
\noindent ограничен  в $X$, так как оценка
$\|V^\lambda_j(u_j)x\|\leq M\|x\|(e^{u_j{\rm Re}\lambda_j}-1)/{\rm
Re}\lambda_j$ влечет
$$
\|W^\lambda_jx\|\leq
\left(c_1^j+M^{n-1}\int\limits_{\Bbb{R}_+^n}\|V^\lambda_j(u_j)\|d\mu(u)\right)\|x\|\leq
$$
$$
\leq\left(c_1^j+\frac{M^n}{{\rm
Re}\lambda_j}\int\limits_{\Bbb{R}_+^n}(e^{u_j{\rm
Re}\lambda_j}-1)d\mu(u)\right)\|x\|=
$$
$$
=\left(c_1^j+\frac{M^n}{{\rm
Re}\lambda_j}(\psi((\mathrm{Re}\lambda_j)
e_j)-c_1^j\mathrm{Re}\lambda_j-\psi(-0))\right)\|x\|
$$
 ($(e_j)$ --- стандартный
базис в $\Bbb{R}^n$). Кроме того, ${\rm Im(}W^\lambda_j)\subseteq
D(A_j)$,  и $W^\lambda_j$ коммутируют со всеми $T_i$, а значит и с
$A_i$. Поэтому, если $\mathrm{Re}\lambda_j<0$ для всех $j$, то из
формулы (2.3)  следует, что
$$
(\psi(\lambda)I-\psi(A))x=\sum_{j=1}^nW^\lambda_j(\lambda_jI-A_j)x,\
x\in D(A).\eqno(2.5)
$$

Пусть теперь $\lambda\in\sigma_a(A)$. Тогда
$\lambda_j\in\sigma(A_j)$, и стало быть (полугруппы $T_j$
ограничены),  $\mathrm{Re}\lambda_j\leq 0$ при всех $j$. Пусть
последовательность  $x_m\in D(A), \|x_m\|=1$ есть совместный
аппроксимативный собственный вектор набора $A$, отвечающий
$\lambda$.

Если $\mathrm{Re}\lambda_j<0$ при всех $j$, то, подставляя в (2.5)
$x_m$ вместо $x$ и полагая $m\to\infty$, выводим, что
$\psi(\lambda)\in\sigma_a(\psi(A))$.

 Теперь предположим, что  $\partial\psi(-0)/\partial s_j\ne\infty$ при всех $j$.
 Заменяя в (2.5) $\lambda$ на $\lambda-\bar\varepsilon$, где
$\bar\varepsilon=(\varepsilon,\dots,\varepsilon)$, $A$ на
$A-\bar\varepsilon I$, а $x$ на $x_m$, получаем
$$
\|(\psi(\lambda-\bar\varepsilon)I-\psi(A-\bar\varepsilon
I))x_m\|\leq
\sum_{j=1}^n\|W^{\lambda-\bar\varepsilon}_j\|\|(\lambda_jI-A_j)x_m\|.
$$
При этом
$$
\|W^{\lambda-\bar\varepsilon}_j\|\leq c_1^j+\frac{M^n}{{\rm
Re}\lambda_j-\varepsilon}(\psi((\mathrm{Re}\lambda_j-\varepsilon)
e_j)-c_1^j(\mathrm{Re}\lambda_j-\varepsilon)-\psi(-0)),
$$
и предел правой части при $\varepsilon\to+0$ существует и конечен
для всех  $j$. Следовательно, существует такая константа $K>0$,
что при всех достаточно малых  $\varepsilon>0$
$$
\|(\psi(\lambda-\bar\varepsilon)I-\psi(A-\bar\varepsilon
I))x_m\|\leq K\sum_{j=1}^n\|(\lambda_jI-A_j)x_m\|. \eqno(2.6)
$$

Далее, $\psi(A-\bar\varepsilon I)x\to \psi(A)x\
(\varepsilon\downarrow+0)$ для любого $x\in D(A)$. Действительно,
при $x\in D(A)$ (оператор $A_j-\varepsilon I$ есть генератор
полугруппы $e^{-\varepsilon u_j}T_j(u_j)$)
$$
\psi(A)x-\psi(A-\bar\varepsilon
I)x=(\psi(-0)-\psi(-\bar\varepsilon))x+
\int\limits_{\Bbb{R}_+^n}(T(u)-I)x (1-e^{-\bar\varepsilon\cdot u})
d\mu(u),
$$
\noindent а потому
$$
\|\psi(A)x-\psi(A-\bar\varepsilon
I)x\|\leq\left(|\psi(-0)-\psi(-\bar\varepsilon)|+
(M^n+1)\int\limits_{\Bbb{R}_+^n} (1-e^{-\bar\varepsilon\cdot u})
d\mu(u)\right)\|x\|,
$$
и осталось применить к правой части  теорему Б. Леви.
Поэтому, полагая в (2.6) $\varepsilon\downarrow +0$, получаем
$$
\|(\psi(\lambda)x_m-\psi(A)x_m\|\leq
K\sum_{j=1}^n\|\lambda_jx_m-A_jx_m\|,
$$
а значит снова $\psi(\lambda)\in\sigma_a(\psi(A))$.

5) Это следует из утверждений 1), 4) и того, что
$\sigma(\psi(A))=\sigma_a(\psi(A))\cup \sigma_R(\psi(A)).$ $\Box$

\bigskip
{\bf 3 Голоморфность полугруппы $g(A)$}

Известно, что условие
$$
\overline{\lim\limits_{t \to +0}}\|I-T(t)\|<2\eqno(3.1)
$$
является достаточным для голоморфности однопараметрической
$C_0$-полу\-группы $T$ в банаховом пространстве $X$. Хотя в общем
случае обратное неверно, это условие и необходимо для
голоморфности  $T$, если $X$ равномерно выпукло (см., например,
[20]). В связи с этим однопараметрическую $C_0$-полугруппу,
удовлетворяющую условию $(3.1)$, будем называть {\it равномерно
голоморфной в} $X$.

Следующая теорема обобщает  утверждения из [9] и [16].

 {\bf 13. Теорема}. {\it Предположим, что полугруппы
$T_j$ удовлетворяют условиям  $\|T_j(t)\|\leq M_j\ (j=1, \dots,
n)$ и
$$
\sum\limits_{j=1}^nC_{j}\overline{\lim\limits_{t \to
+0}}\|I-T_j(t)\|<2,
$$
где $C_j=\prod_{k=1}^{j-1}M_k$ при $j>1,\ C_1=1$.  Тогда для любой
функции $\psi$ из ${\cal T}_n$ оператор $\psi(A)$ является
генератором равномерно голоморфной полугруппы.}

Доказательство. Заменяя, если это необходимо, $\psi$ на
$\psi-\psi(-0)$, можем считать, что $c_0(=\psi(-0))=0$. Положим
$b_j=\overline{\lim}_{t \to +0}\|I-T_j(t)\|$ и выберем
$\varepsilon > 0$ таким, что $\sum_{j=1}^n C_{j}b_j+\varepsilon
<2$. Найдется такое $\delta >0$, что $\|I-T_j(t)\|<b_j+\varepsilon
/(nC_{j})$ при всех $j=1,\dots ,n; t\in [0;\delta)$.

Далее,  из тождества ($T_0(t):=I$)
$$
T(u)-I=\sum\limits_{j=1}^n\left(\prod_{k=0}^{j-1}
T_{k}(u_k)\right)(T_j(u_j)-I)
$$
следует, что

$$
\|I-T(u)\|\leq \sum\limits_{j=1}^nC_{j}\|I-T_j(u_j)\|,\eqno(3.2)
$$
а потому при $u\in[0;\delta)^n$ справедливо неравенство
$$\|I-T(u)\|\leq \sum\limits_{j=1}^n C_{j}b_j+\varepsilon.
$$
 Следовательно,
если $x\in X, \|x\|= 1$, то
$$
\|(I-g_t(A))x\|\leq
\int\limits_{\Bbb{R}_+^n}\|I-T(u)\|d\nu_t(u)\|x\|=
$$
$$
=\int\limits_{[0;\delta)^n}\|I-T(u)\|d\nu_t(u)+\int\limits_{\Bbb{R}_+^n\setminus
[0;\delta)^n}\|I-T(u)\|d\nu_t(u)\leq
$$
$$
\leq\sum_{j=1}^n
C_{j}b_j+\varepsilon+\int\limits_{\Bbb{R}_+^n\setminus
[0;\delta)^n}\|I-T(u)\|d\nu_t(u).
$$
То есть
$$
\|I-g_t(A)\|\leq \sum_{j=1}^n
C_{j}b_j+\varepsilon+\int\limits_{\Bbb{R}_+^n\setminus
[0;\delta)^n}\|I-T(u)\|d\nu_t(u).\eqno(3.3)
$$

 Заметим теперь, что направленность мер $\nu_t$ узко
сходится к мере Дирака $\varepsilon_0$ при $t\to +0$. В самом
деле, преобразование Лапласа ${\cal L}\nu_t(s)=e^{t\psi(s)}$
непрерывно в точке $s=0$ и ${\cal L}\nu_t(s)\to 1={\cal
L}\varepsilon_0(s)$ при $t\to +0$. Поэтому узкая сходимость
вытекает из теоремы непрерывности для многомерного преобразования
Лапласа (см., например, [19, глава~IX, \S\,5, теорема 3 с]). Но
так как полугруппы $T_j$ удовлетворяют условию
$\overline{\lim\limits_{t \to +0}}\|I-T_j(t)\|<2$, они голоморфны
(см., например, [20, следствие 2.5.7]), а потому становятся
непрерывными в топологии нормы. Значит,  ограниченная
 функция $u\mapsto\|I-T(u)\|$ непрерывна на $\Bbb{R}_+^n\setminus
[0;\delta)^n$. Следовательно, переходя в (3.3) к верхнему пределу
при $t\to +0$, получим
$$\overline{\lim\limits_{t \to
+0}}\|I-g_t(A)\|\leq\sum_{j=1}^n C_{j}b_j+\varepsilon<2.
$$
В силу отмеченного выше критерия отсюда следует равномерная
голоморфность полугруппы $g(A)$, что и требовалось доказать.$\Box$

{\bf 14. Следствие}. {\it Пусть пространство $X$ равномерно
выпукло, $T_1$ -- голоморфная полугруппа в $X$, а операторы
$A_2,\dots,A_n$ ограничены (если $n>1$). Тогда для любой функции
$\psi$ из ${\cal T}_n$ оператор $\psi(A)$ является генератором
голоморфной полугруппы.}

\medskip
Доказательство. Условие теоремы выполнено, поскольку
$\overline{\lim}_{t \to +0}\|I-T_1(t)\|<2$ (см., например, [20,
следствие 2.5.8]), и при $j>1$ справедливы равенства $\lim_{t \to
+0}\|I-T_j(t)\|=0$. $\Box$

\bigskip
{\bf 15. Следствие}. {\it Пусть пространство $X$ равномерно
выпукло. Если $T$ -- однопараметрическая голоморфная полугруппа
сжатий в $X$ с генератором $A$, то  для любой функции $\psi$ из
${\cal T}_1$ оператор $\psi(A)$ является генератором голоморфной
полугруппы сжатий.}

Следствие 15 дает для случая равномерно выпуклых пространств
положительный ответ на один вопрос из [16].

\bigskip
{\bf 4 Неравенство моментов}

Имеет место следующее неравенство, обобщающие результаты из [21] и
[16].

{\bf 16. Теорема}. {\it Если  операторы $A_j\in {\rm Gen}(X)$
порождают $C_0$-полугруппы $T_j$ соответственно, причем
$\|T_j(t)\|\leq M\ (t\in \Bbb{R}_+; j=1, \dots, n; M\geq 1)$,  то
для любой функции $\psi\in {\cal T}_n$ и любого $x\in D(A),
\|x\|=1$ справедливо неравенство}
$$
\|\psi(A)x\|\leq -nK_MM^{n-1}\psi\left(-\frac{1}{n}\|A_1x\|,\dots,
-\frac{1}{n}\|A_nx\|\right),
$$
\textit{где} $K_M=(M+1)/(1-e^{-(M+1)/M})$.

Доказательство. Легко проверить, что теорема верна для линейных
функций Бернштейна $\psi(s)=c_0+c_1\cdot s$, поскольку $K_M>1$.
Кроме того, в силу формулы (1.1) $\psi(s)=c_0+c_1\cdot s+
\psi_0(s)$, где $\psi_0\in {\cal T}_n$ и не содержит линейной
части. Следовательно, можно считать, что в (1.1) $c_0=c_1=0$.
Формулы (1.2) и (1.1) показывают, что достаточно доказать
неравенство
$$
\|(T(u)-I)x\|\leq nK_M M^{n-1}\left(1-e^{-\|Ax\|\cdot
u/n}\right),\quad u>0,\eqno(4.1)
$$
где положено $\|Ax\|\cdot u=\sum_{j=1}^n \|A_jx\|u_j$.

Докажем его сначала  при $n=1$, полагая для краткости $A_1=A,
T_1(u_1)=T(u)$.
 Обозначим через $t(r)$ функцию, обратную возрастающей функции
$r(t)=t/(1-e^{-t}),\ r(0)=1$. Для фиксированного $r\geq 1$
возможны два случая.

1) $\|Ax\|u\leq t(r)$. Тогда $r(\|Ax\|u)\leq r(t(r))=r$, т. е.
$\|Ax\|u \leq r\cdot\left(1-e^{-\|Ax\|u}\right)$, а потому
$$
\|(T(u)-I)x\|=\left\|\int\limits_0^1\frac{d}{ds}T(us)xds\right\|=
$$
$$
\left\|\int\limits_0^1T(us)Axuds\right\|\leq M\|Ax\|u\leq
Mr\cdot\left(1-e^{-\|Ax\|u}\right).
$$

2) $\|Ax\|u>t(r)$. Тогда
$$
\|(T(u)-I)x\|\leq
M+1\leq\frac{M+1}{1-e^{-t(r)}}\left(1-e^{-\|Ax\|u}\right).
$$
В любом случае справедливо неравенство
$$
\|(T(u)-I)x\|\leq C(r)\left(1-e^{-\|Ax\|u}\right),
$$
где $C(r)=M\max\{r;(M+1)/M\left(1-e^{-t(r)}\right)\}$.
Для минимизации $C(r)$ заметим, что функция $t(r)$ возрастает от 0
до $+\infty$ при $1\leq r<+\infty$. Поэтому уравнение
$r=(M+1)/M\left(1-e^{-t(r)}\right)$, т. е. $t(r)=(M+1)/M$, имеет
единственное решение $r_0=r((M+1)/M)$. Если $r<r_0$, то в силу
отмеченной монотонности
$$
\frac{M+1}{M\left(1-e^{-t(r)}\right)}>
\frac{M+1}{M\left(1-e^{-t(r_0)}\right)}=\frac{M+1}{M\left(1-e^{-(M+1)/M)}\right)}=r_0.
$$
Таким образом,
$$
\min\{C(r):r\geq 1\}=Mr_0=(M+1)/\left(1-e^{-(M+1)/M}\right),
$$
что и доказывает  (4.1) при $n=1$.

Перейдем к общему  случаю. Используя неравенство (3.2),
 заключаем, что
$$
\|I-T(u)\|\leq \sum\limits_{j=1}^nM^{j-1}\|I-T_j(u_j)\|\leq
M^{n-1}\sum\limits_{j=1}^n\|I-T_j(u_j)\|.
$$

Применяя здесь к каждому слагаемому в правой части доказанный выше
частный случай  неравенства (4.1) и воспользовавшись неравенством
Коши, получаем
$$
\|I-T(u)\|\leq
K_MM^{n-1}\sum\limits_{j=1}^n\left(1-e^{-\|A_jx\|u_j}\right)\leq
nK_MM^{n-1}\left(1-e^{-\frac{1}{n}\sum\limits_{j=1}^n\|A_jx\|u_j}\right),
$$
и теорема доказана. $\Box$

{\bf 17. Следствие}. {\it Для любой функции $\psi\in {\cal T}_n$ и
любого вектора $x\in D(A), x\ne 0$ справедливо неравенство}
$$
\|\psi(A)x\|\leq
-nK_MM^{n-1}\psi\left(-\frac{\|A_1x\|}{n\|x\|},\dots,
-\frac{\|A_nx\|}{n\|x\|}\right)\|x\|.
$$

Для формулировки другого следствия выделим некоторый класс
функциональных пространств. Будем говорить, что комплексное
банахово пространство $X$ принадлежит классу ${\cal C}$, если оно
обладает следующими свойствами:

а) $X$ есть пространство функций на $\mathbb{R}^n$, содержащее
пространство основных функций ${\cal D}(\mathbb{R}^n)$  в качестве
векторного подпространства.

б) Если последовательность функций $x_k\in {\cal
D}(\mathbb{R}^n)$, носители которых содержатся в общем компакте из
$\mathbb{R}^n$, сходится равномерно на $\mathbb{R}^n$ к $0$, то
$x_k$ сходится к $0$ и в $X$.

в) Если последовательность функций $x_k\in X$ сходится в $X$ к
$0$, то она содержит такую подпоследовательность $x_{k_m}$, что
$x_{k_m}(0)\to  0$.

г) Для любой функции $x\in X$ и любого $u\in \mathbb{R}^n_+$ сдвиг
$x_u$, где $x_u(s)=x(u+s)$, также принадлежит $X$, причем
операторы сдвига $T(u)x=x_u$  ограничены в $X$ и образуют
$n$-параметрическую ограниченную $C_0$-полугруппу.

Пространства  $C_0(\mathbb{R}^n)$ непрерывных функций на
$\mathbb{R}^n$, исчезающих на бесконечности, и $UCB(\mathbb{R}^n)$
равномерно непрерывных ограниченных  функций на $\mathbb{R}^n$,
наделенные $\sup$-нормой, принадлежат классу ${\cal C}$.

{\bf 18. Следствие}. {\it  Если функция $\psi\in {\cal T}_n$
ограничена на $(-\infty;0)^n$, то для любого банахова пространства
$X$ оператор $\psi(A)$ ограничен при всех попарно коммутирующих
$A_j\in {\rm Gen}(X)$ (причем при фиксированной $\psi$ его
константа ограниченности зависит лишь от $M$). Обратно, если для
некоторого банахова пространства $X$ класса ${\cal C}$ оператор
$\psi(A)$, где $A=(\partial_1,\dots,\partial_n)$ --- набор частных
генераторов полугруппы сдвигов в $X$, ограничен, то функция $\psi$
ограничена на $(-\infty;0)^n$.}

Доказательство. Первое утверждение с очевидностью следует из
неравенства моментов. Для доказательства обратного заметим, что
частными генераторами полугруппы сдвигов $T(u)$ в пространстве $X$
служат операторы $A_j=\partial_j$  взятия частных производных,
причем их область определения содержит пространство ${\cal
D}(\mathbb{R}^n)$ основных функций (с учетом свойств а), б) это
доказывается так же, как и в случае  $X=C_0(\mathbb{R}^n)$).
Предположим, что оператор $\psi(A)$ ограничен, и докажем, что мера
$\mu$ на $\mathbb{R}_+^n$ (см. (1.1)) конечна. Допустим, что это
не так. Поскольку в силу [3, лемма 3.1]
$\mu(\mathbb{R}_+^n\setminus[0,\delta)^n)<\infty$ при всех
$\delta>0$, из этого допущения следует, что
$\mu(V_k\cap\mathbb{R}_+^n)\to\infty\ (k\to\infty)$, где
$V_k=\{s\in \mathbb{R}^n:1/k\leq |s|\leq 1\}$ --- полый шар. Для
каждого натурального $k$ выберем функцию  $x_k\in{\cal
D}(\mathbb{R}^n)$ так, что $0\leq x_k\leq
1/\mu(V_k\cap\mathbb{R}_+^n)$, носители ${\rm supp}(x_k)$
содержатся в шаре $B[0,3]\subset \mathbb{R}^n$, причем $x_k=0$ в
некоторой окрестности нуля, и сужение
$x_k|V_k=1/\mu(V_k\cap\mathbb{R}_+^n)$. Тогда $x_k\to 0$ в $X$ в
силу б). Далее, из (1.2) следует, что
$$
\psi(A)x_k(s)=c_0x_k(s)+c_1\cdot Ax_k(s)+\int\limits_{\Bbb{R}_+^n}
(x_k(s+u)-x_k(s))d\mu(u).
$$
Так как $x_k(0)=\partial_jx_k(0)=0\ (j=1,\dots,n)$, имеем
$$
\psi(A)x_k(0)=\int\limits_{\Bbb{R}_+^n} x_k(u)d\mu(u)\geq
\frac{1}{\mu(V_k\cap\mathbb{R}_+^n)}\mu(V_k\cap\mathbb{R}_+^n)=1,
$$
а это противоречит тому, что $\psi(A)x_k\to 0$ в $X$ в силу в).
Таким образом, мера $\mu$  конечна. Но тогда оператор (см. (1.2))
$$
c_1\cdot A=\psi(A)-c_0-\int\limits_{\Bbb{R}_+^n} (T(u)-I)d\mu(u)
$$
ограничен в $X$ вместе с $\psi(A)$. Покажем, что это возможно лишь
при $c_1=0$. Для каждого натурального $k$ рассмотрим функцию
одного переменного
$\varphi_k(t)=k^{-1}\exp(-\frac{k^-2}{k^{-2}-t^2})$ при $|t|\leq
k^{-1}$, $\varphi_k(t)=0$  при $|t|>k^{-1}$. Выберем также такую
функцию $\varphi\in {\cal D}(\mathbb{R}^n)$, что $0\leq\varphi\leq
1,\ {\rm supp}(\varphi)\subset B[0,2]\subset \mathbb{R}^n,\
\varphi|B[0,1]=1$, и для любого $l=1,\dots, n$ положим
$x_k^l(s)=\varphi_k(s_l-k^{-2})\varphi(s)$. Тогда $x_k^l\in{\cal
D}(\mathbb{R}^n),\ {\rm supp}(x_k^l)\subset B[0,2]$, и $x_k^l$
равномерно  стремится к $0$ при $k\to\infty$, а потому $x_k^l\to
0$ в $X$. Поскольку оператор $c_1\cdot
A=\sum_{j=1}^nc_1^j\partial_j$ ограничен в $X$, то и
$\sum_{j=1}^nc_1^j\partial_jx_k^l\to 0$ в $X$ при $k\to\infty$.
Переходя к подпоследовательности, получаем, что
$\sum_{j=1}^nc_1^j\partial_jx_k^l(0)\to 0$ при $k\to\infty$. С
другой стороны, непосредственно проверяется, что
$\partial_lx_k^l(0)\to\infty$ при $k\to\infty$, и
$\partial_jx_k^l(0)=0$ при $j\ne l$. Поэтому $c_1^l=0$.

Следовательно, $c_1=0$ и функция
$$
\psi(s)=c_0+\int\limits_{\Bbb{R}_+^n} (e^{us}-1)d\mu(u)
$$
ограничена на $(-\infty;0)^n$. $\Box$

{\bf 19. Следствие}. {\it Если последовательность функций
$\psi_k\in {\cal T}_n$ сходится к нулю поточечно на
$(-\infty;0]^n$, то $\psi_k(A)x\to 0$ при всех $x\in D(A)$.}

\centerline{Аннотация}

\vspace{5mm}
 Развивается  многомерное
функциональное исчисление генераторов полугрупп, основанное  на
классе функций Бернштейна нескольких переменных. Устанавливаются
теоремы об отображении спектров, дается условие голоморфности
полугрупп, порождаемых операторами, возникающими в исчислении,  а
также  доказывается неравенство моментов для этих операторов.

\hspace{5mm}

{\it Ключевые слова: многопараметрическая полугруппа операторов,
голоморфная полугруппа, генератор полугруппы, функциональное
исчисление, функция Бернштейна, теорема об отображении спектра}

\begin{center}
{\bf ON SOME PROPERTIES  OF MULTIDIMENSIONAL BOCHNER-PHILLIPS
FUNCTIONAL CALCULUS \\
 A. R. Mirotin}\\
\end{center}

\hspace{5mm}

\centerline{Abstract}

The multidimensional functional calculus of semigroup generators,
based on the class of Bernstein functions in several variables is
developed, the spectral mapping
 theorems for joint  spectra have been stated,
 the condition for holomorphy  of semigroups, generated by
operators which arises in the calculus  is given,  and the
 moment inequality for such operators is proved.

\hspace{5mm}

{\it Key words and phrases: multiparameter semigroup of operators,
holomorphic semigroup, semigroup generator, functional calculus,
Bernstein function, spectral mapping theorem}

\end{document}